\documentclass{amsart}
\usepackage{amsmath,amssymb,amsbsy,amsfonts,amsthm,latexsym,
                       amsopn,amstext,amsxtra,euscript,amscd}
\usepackage{array}
\usepackage{ifthen}
\usepackage{url}
\theoremstyle{plain}
\newtheorem{theorem}{Theorem}[section]
\newtheorem{prop}[theorem]{Proposition}
\newtheorem{lemma}{Lemma}[section]
\newtheorem{corol}{Corollary}[theorem]
\newtheorem{conj}{Conjecture}[section]
\theoremstyle{definition}

\theoremstyle{remark}

\numberwithin{equation}{section}

\begin{document}
\title[On the conjecture of Je\'{s}manowicz]
{On the conjecture of Je\'{s}manowicz}

\author{G\"{o}khan Soydan, Musa Dem\.irc\.i, Ismail Naci Cangul, and Alain Togb\'{e} }

\address{{\bf G\"{o}khan Soydan} \\
Department of Mathematics \\
Uluda\u{g} University\\
 16059 Bursa, Turkey}
\email{gsoydan@uludag.edu.tr }

\address{{\bf Musa Dem\.irc\.i}\\
Department of Mathematics \\
Uluda\u{g} University\\
 16059 Bursa, Turkey}
\email{mdemirci@uludag.edu.tr}

\address{{\bf \.{I}smail Nac\.i Cang\"{u}l}\\
Department of Mathematics \\
Uluda\u{g} University\\
 16059 Bursa, Turkey}
\email{cangul@uludag.edu.tr}

\address{ {\bf Alain Togb\'e}\\
Department of Mathematics, Statistics, and Computer Science, Purdue University Northwest\\
 1401 S, U.S. 421, Westville IN 46391, USA}
\email{atogbe@pnw.edu}

\newcommand{\acr}{\newline\indent}

\thanks{}

\subjclass[2010]{11D61}
\keywords{Je\'{s}manowicz' conjecture, exponential Diophantine equation}

\begin{abstract}
We give a survey on some results covering the last 60 years concerning Je\'{s}manowicz' conjecture. Moreover, we conclude the survey with a new result by showing that the special Diophantine equation  $$(20k)^x+(99k)^y=(101k)^z$$ has no solution other than $(x,y,z)=(2,2,2)$.
\end{abstract}

\date{\today}
\maketitle

%######################################################%
\section{The life and work of Leon Je\'{s}manowicz}\label{sec:1}
%######################################################%

Professor Leon Je\'{s}manowicz was born in $27.04.1914$ in Druja in the Vilnius region. He is the son of a postal clerk Anatole Je\'{s}manowicz and Irene. In 1920, his parents moved to Vilnius. In the same year his father died and his mother moved first to {\L}odzi and then to Grodno. Leon Je\'{s}manowicz graduated from high school in the area of humanities in 1932.\\

During 1933-1937, he studied mathematics at the Stefan Batory University in Vilnius. The first two years, he also studied drawing in the Faculty of Fine Arts in addition to Mathematics. After graduation, he received a scholarship from the National Culture Fund and started working as a junior assistant at the same University. Until the outbreak of world war II, he worked as an assistant in the Mathematics Department under the guidance of Professor Antoni Zygmund.\\ 

In 1939, he prepared his doctoral thesis, but he had no chance to defend it because of the outbreak of war. During the German occupation, he conducted clandestine study groups. In March 1945, he returned back to Lublin with his family where he got a post in the Department of Mathematics at Maria Curie-Sk{\l}odowska University. He finally defended his doctoral dissertation on the uniqueness of Schl\"{o}milcha series in July 1945. His supervisor was Prof. Juliusz Rudnicki. From 1.10.1945, he worked as a senior assistant at the same University. In 1946, he permanently settled in Torun, where he was associated with the newly established Mathematics Department at the Nicolaus Copernicus University until his death in 1989.\\ 

He became an assistant professor in 1949, an associate professor in 1954, and in 1964 the State Council awarded him the title of professor. For many years, he held positions in Faculty of Mathematics, Physics and Chemistry including vice-dean and dean. He also performed many social duties. He was a long-standing chairman of the Torun branch of the Polish Mathematical Society and member of the board.\\ 

He studied the theory of Sumawalno\'{s}ci ranks (see \cite{Jes1},  \cite{Jes2}, \cite{Jes4}, \cite{Jes5}, \cite{Jes7}), the theory of abelian groups (see \cite{Jes6}), and the solutions of Diophantine equations (see \cite{Je}).\\ 

Professor Je\'{s}manowicz had many varied interests apart from Mathematics. He was also interested in literature, history, theatre, and caricature. In his drawings, he represented hundreds of people and groups with whom he had contact during his life (see \cite{BoKa} and \cite{Int}).

%##############################################%
\section{The conjecture of Je\'{s}manowicz} \label{sec:2}
%###############################################%

In 1955/56, Sierpin\'{s}ki \cite{Si} showed that the equation $3^{x}+4^{y}=5^{z}$ has $x=y=z=2$ as its only solution in positive integers. The integers $3,4,5$ constitute a triple of Pythagorean numbers, integral solutions of the equation $a^{2}+b^{2}=c^{2}$. In the same year, Je\'{s}manowicz \cite{Je} proved that Sierpin\'{s}ki's result holds also for the following set of Pythagorean numbers:
\begin{equation}\label{eq:2.1}
2n+1, 2n(n+1), 2n(n+1)+1,
\end{equation}
with $n = $ 2, 3, 4, 5 (note that $n=1$ is the Sierpin\'{s}ki's case).\\

Let $U,V,W$ be fixed positive integers. Consider the Diophantine equation
\begin{equation}\label{eq:2.2}
U^{x}+V^{y}=W^{z}.
\end{equation}
Je\'{s}manowicz proposed the following problem:
\begin{conj}[\textbf{Je\'{s}manowicz' Conjecture}]
Assume $U^{2}+V^{2}=W^{2}$. Then, equation \eqref{eq:2.2} has no positive solution $(x,y,z)$ other than $x=y=z=2.$ 
\end{conj}
It is well-known that the numbers
\begin{equation}\label{eq:2.3}
U=p^{2}-q^{2}, V=2pq, W=p^{2}+q^{2} 
\end{equation}
form all solutions for 
\begin{equation}\label{eq:2.4}
U^{2}+V^{2}=W^{2},
\end{equation}
where $(p,q)=1$, $p>q$, $p$ and $q$ have opposite parity.

The first results on Je\'{s}manowicz' conjecture were obtained by Ko \cite{CKo6} in $1958$ without calling it Je\'{s}manowicz' conjecture. Ko proved the conjecture if\\

$(i)$ $n\equiv$ 1, 4, 5, 9, 10 $\pmod{12}$;\\

$(ii)$ $n$ is odd and there exist a prime $p$ and a positive integer $s$ such that $2n+1=p^{s}$,\\
%and\\

$(iii)$ a prime $p\equiv 3\pmod{4}$ and an integer $n$ which is the sum of two squares, such that $2n+1\equiv 0\pmod{p}$.\\
 
In the same year, Ko \cite{CKo1} proved the conjecture for several cases by giving the following three theorems:

\begin{theorem}
$(i)$ Je\'{s}manowicz' conjecture is true for $n\equiv$ 3, 7, 11 $\pmod{12}$.\\

$(ii)$ If there exists  prime $p\equiv 3\pmod{4}$ such that $2n+1\equiv 0\pmod{p}$, then the conjecture holds.\\

$(iii)$ If there exists a prime $p\equiv 5\pmod{8}$ such that $2n+1\equiv 0\pmod{p}$, then the conjecture holds for all integers $a, \ b, \ c.$
\end{theorem}

\begin{theorem}
Je\'{s}manowicz' conjecture  holds for $n\equiv 2\pmod{5}$, $n\equiv 3\pmod{7}$, $n\equiv 4\pmod{9}$, $n\equiv 5\pmod{11}$, $n\equiv 6\pmod{13}$, and $n\equiv 7\pmod{15}$.
\end{theorem}

\begin{theorem}
When $n<96$, Je\'{s}manowicz' conjecture  holds for all integers $a,b,c$. 
\end{theorem}

A year later, Ko \cite{CKo2} proved another theorem covering same special cases of the conjecture:

\begin{theorem}
$\textbf{(I)}$ In \eqref{eq:2.3}, if the numbers $p=2n$ and $q$ contain no prime factor congruent to $1$ modulo $4$, $2n>q>0$, $(2n,q)=1$, and if one of the following conditions holds, then the conjecture is true:\\

$(i) $ $n\equiv 2\pmod{4}$, $q\equiv 3\pmod{8}$,\\

$(ii)$ $n\equiv 2\pmod{4}$, $q\equiv 5\pmod{8}$, $2n+q$ has a prime factor congruent to $3$ modulo $4$,\\

$(iii)$ $n\equiv 0\pmod{4}$, $q\equiv 3,5\pmod{8}$.\\

$\textbf{(II)}$ If $p=3n$ and $q=2m$ have no prime factors congruent to $1$ modulo $4$, $(3n,2m)=1$, $2\sqrt{2}m>3n>2m>0$ or $3n>8m>0$ and if one of the following conditions holds, then the conjecture is true:\\

$(i)$ $m\equiv 2\pmod{4}$, $n\equiv 1\pmod{8}$,\\

$(ii)$ $m\equiv 2\pmod{4}$, $n\equiv 7\pmod{8}$, $3n+2m$ has a prime factor congruent to $3$ modulo $4$,\\

$(iii)$ $m\equiv 0\pmod{4}$, $n\equiv 1,7\pmod{8}$.
\end{theorem}
The same year, Lu \cite{Lu1} considered the case $p=2n$, $q=1$ of \eqref{eq:2.3} and showed that the conjecture still holds in this case.\\ 

In $1961$, J\'{o}zefiak \cite{Joz} confirmed the conjecture for a class of Pythagorean numbers: If $U=2^{2r}p^{2s}-1$, $V=2^{r+1}p^{s}$, $W=2^{2r}p^{2s}+1$, where $r,s\in\mathbb{N}$, $\mathbb{N}=\{0,1,2,...\}$ denotes the set of natural numbers and $p$ is a prime number, then $x=y=z=2$ is the only integral solution of the equation \eqref{eq:2.2}.\\ 

In $1962$, Podyspanin \cite{Pody} proved that equation \eqref{eq:2.2} for $U, \ V, \ W$ given as in \eqref{eq:2.1} has no solution in positive integers $n, \ x, \ y, \ z$ when $4\nmid n$ or when $2n+1$ has a factor not congruent to $1$ modulo $8$. He obtained infinitely many primitive Pythagorean triples $U, \ V, \ W$ such that equation \eqref{eq:2.2} has no solutions $x, \ y, \ z$ in natural numbers except $x=y=z=2$, which is a result already proved by  Lu \cite{Lu1}.\\ 

In $1963$, Ko \cite{CKo3} considered the values of $n$ modulo $240$ and proved the conjecture for the following values of $n$ in \eqref{eq:2.1}:\\

$(i)$ $n \not\equiv$ 0, 24, 48, 80, 96, 104, 120, 128, 144, 176, 200, 224  $\pmod{240}$,\\

$(ii)$ $n \equiv$ 0, 24, 80, 104, 120, 144, 200, 224  $\pmod{240}$ and there exists a prime  $p \not\equiv 1 \pmod{16}$ and $p \mid (2n+1)$,\\

$(iii)$ $n \equiv$ 48, 96, 128, 176 $\pmod{240}$ and there exists a prime  $p \not\equiv 1 \pmod{32}$ and $p \mid (2n+1)$,\\

$(iv)$ $n \leq 1000$ and $n\neq$ 96, 120, 128, 144, 200, 224, 288, 320, 336, 384, 440, 464, 564, 576, 600, 608, 624, 680, 704, 744, 800, 914, 960.\\

In $1964$, Ko and Sun \cite{CKo4} used similar arguments and techniques when the numbers are taken modulo $128$ and generalized the above $(iv)$ to show that the Je\'{s}manowicz'  conjecture holds for all $n\leq 1000$. They also gave some special classes of the values of $n$ for which the conjecture is true.\\

The same year, in the next issue of the same journal, Ko \cite{CKo5} improved the bound to $n<6144$ and also gave some classes of the values of $n$ for which the conjecture holds.\\

In $1965$, Dem'janenko \cite{Dem} proved that equation \eqref{eq:2.2} has no solution for $W=V+1$ and for $U=m^{2}-1$, $V=2m$, $W=m^{2}+1$, where $m \in \mathbb{N^{+}}$ using the results in \cite{CKo4} and \cite{Pody}.\\

After a long silence, in $1984$, Grytczuk and Grelak \cite{GrGr} proved the conjecture for two different cases:\\

$(i)$ For $p=2a$ and $q=1$ in \eqref{eq:2.4};\\

$(ii)$ For $U=2^{\alpha}-1$, $V=2(2^{\alpha}+1)$ and $W=3.2^{\alpha}+1$ with $2^{\alpha}+1$ is prime. Actually, the first case was also proved by Lu \cite{Lu1} in $1959$.\\

In $1993$, Deng \cite{De2} proved the conjecture for some even $V$, where some congruence conditions are satisfied for the divisors of $V$.\\

The same year, Takakuwa and Asaeda \cite{Takas} tried to obtain a generalization of the triples satisfying \eqref{eq:2.4} to triples satisfying \eqref{eq:2.2}. They considered the case where $p=2a$ and $q\equiv 3 \pmod{4}$ is a prime and studied the Diophantine equation
\begin{equation*}
(4a^{2}-q^{2})^{l}+(4aq)^{m}=(4a^{2}+q^{2})^{n},
\end{equation*}
where $a,q\in \mathbb{N^{+}}$, with $(a,q)=1$, $2a>q$. 
One must notice that this equation is a generalization of the equation considered by Grytczuk and Grelak \cite{GrGr} with $q=1$. They proved the conjecture for this triple by considering the different parities of $a$ and also given some conditions on the divisors of $a$ for the conjecture to be true.  Particularly, they proved it when $a$ is odd, $q$ is an odd prime, $m\neq 1$ and $q\equiv 3 \pmod{4}$; $a$ is even, $q\equiv 3 \pmod{8}$ is an odd prime, $2a+q$ is prime and $2a-q$ is prime or $1$; and finally $a$ is odd, $q\equiv 3 \pmod{4}$ is an odd prime, and when a prime divisor $p$ of $a$ satisfies the conditions $p\equiv 1 \pmod{4}$ and $(\frac{q}{p})=-1$. They claimed that their method of proof also covers some other specific cases under some given conditions on $a$ and $q$.\\

In $1993$, Terai \cite{Ter1} proposed an analogue of Je\'{s}manowicz' conjecture:
\begin{conj}\label{333}
If $a^{2}+b^{2}=c^{2}$ with $(a,b,c)=1$ and $a$ is even, then the Diophantine equation
\begin{equation}\label{444}
x^{2}+b^{m}=c^{n}
\end{equation}
has the only positive integer solution $(x,m,n)=(a,2,2)$. 
\end{conj}

In the same paper, Terai proved his conjecture for primes $b$ and $c$ with\\

$i)$ $b^{2}+1=2c$,\\ 

$ii)$ $d=1$ or even if $b\equiv 1 \pmod{4}$, where $d$ is the order of a prime divisor of $[c]$ in the ideal class group $\mathbb{Q}(\sqrt{-b})$. Further, he proved that his conjecture holds when $b^{2}+1=2c$, $b<20$, $c<200$.\\

Later in $1998$, Cao and Dong \cite{CaDo2} proved that if $(i)$ $b$ is a prime power, $c\equiv 5\pmod{8}$, or $(ii)$ $c\equiv 5\pmod{8}$ is a prime power, then Conjecture \ref{333} holds.\\

In $2001$, Cao and Dong \cite{CaDo} proved that if 
$$a=V_{r}, \ b=U_{r}, \ c=m^{2}+1, \ b>8.10^{6} \ \mbox{and} \  b \equiv 3 \pmod{4},$$
then the Diophantine equation \eqref{444} has only the positive integer solution $(x,m,n)=(a,2,r)$, where $b$ is an odd prime, $m,r\in \mathbb{N}$, $2\mid m$, $2\nmid r$, $r>1$ and 
$$(m+i)^{r}=V_{r}+i \cdot U_{r}, \ i=\sqrt{-1}.$$

There are several papers dealing with Terai conjecture, see e.g. \cite{CaDoLi}, \cite{CaDo}, \cite{De3}, \cite{HuLe}, \cite{Le1}, \cite{Le2}, \cite{Ter2}, \cite{TerTa} and \cite{YuWa}. A brief discussion on the results concerning Terai conjecture will be given in Section $4$.\\

In $1994$, Chen \cite{Chen} studied the numbers in \eqref{eq:2.1} and prove a special case of the Je\'{s}manowicz' conjecture where $2n+1\equiv 0 \pmod{p}$ with $p\equiv 3 \pmod{4}$ is prime, and where $n\equiv 1 \pmod{3}$. Another specific result proving Je\'{s}manowicz' conjecture is obtained for equation \eqref{eq:2.2} in the particular case where $q=3$ and $p\leq6000$ is even by Guo and Le \cite{GuLe}. The same year, Le \cite{M2} proved the conjecture for the case where $2\parallel pq$, $W=r^{n}$ with $r$ is an odd prime, $n \in \mathbb{N^{+}}$. In \cite{WaDe}, Wang and Deng proved the conjecture for the pythagorean triple \eqref{eq:2.3} where $p$ and $q$ have no prime divisor $r$ with  $r\equiv 1 \pmod{4}$ and $p$ and $q$ satisfy certain congruence conditions.\\

Later in $1996$, Le \cite{M3} using Baker's method for the first time, proved another particular case where  $2\parallel p$, $p\geq 81q$ and $q\equiv 3 \pmod{4}$.\\

Just a short time after Le's paper, Takakuwa \cite{Ta} eliminated the condition $p\geq 81q$ when $q=$ 3, 7, 11, and 15.\\

In $1998$, Chen \cite{Chen2} solved some cases of the conjecture where $W=V+1$ and $(U,V,W)=(m^{2}-1,2mn,m^{2}+1)$.\\ 

Deng and Cohen \cite{DeCo2} proved the conjecture when one of $p$ and $q$ in \eqref{eq:2.3} has no prime factor congruent to $1$ modulo $4$ and also when certain congruence relations on $p$ and $q$ are satisfied.\\

For some time after $1999$, most of the papers on the subject were dealing with several variations of Terai's conjecture. In $2009$, Le \cite{M4} proved the  Je\'{s}manowicz' conjecture under the restriction
\begin{equation*}
\gcd\left(\frac{U^{d}+(-1)^{s}}{W},W\right)=1,
\end{equation*}
where $d$ is the least integer such that 
$W\mid (V^{d}+(-1)^{s}),$
for $s\in \{0,1\}$ chosen to minimize $d$. Under this restriction, he showed that the conjecture holds for all $W>4.10^{9}$.\\

The same year, Miyazaki \cite{Mi1} obtained two new results, similar to those of Deng and Cohen. The first result settled the Je\'{s}manowicz' conjecture when $p^{2}-q^{2}$ has no prime factor congruent to $1$ modulo $4$, $p-q$ has a prime factor congruent to $3$ modulo $8$ and finally $p\not\equiv 1\pmod{4}$. The second dealt with the conjecture when $p\equiv 4\pmod{8}$ and $q\equiv 7\pmod{16}$ or  $p\equiv 7\pmod{16}$ and $q\equiv 4\pmod{8}$, where $p$ and $q$ are as in \eqref{eq:2.3}. Both results also included the case where $V$ is divisible by $8$, unlike the most attempts earlier including only the case $4\mid V$. He gave the following result which is used in the proof of the main theorem and later in \cite{Mi2}, \cite{Mi3}.
\begin{prop}
If $x,y,z$ in equation \eqref{eq:2.2} are all even, then $x/2$, $y/2$ and $z/2$ are all odd.
\end{prop}

In $2010$, Hu and Yuan \cite{HY} proved the conjecture for the triple 
$(U,V,W)=(2n+1,2n(n+1),2n(n+1)+1),$
by the use of the BHV theorem \cite{BHV}, which states the existence conditions of primitive prime divisors of Lucas and Lehmer numbers.\\

In $2011$, Miyazaki \cite{Mi2} considered the case $q>1$ as Lu \cite{Lu1} proved the conjecture for $q=1$. Also he considered the case where the $2$-adic valuation of $p,q$ is $>1$ and took $p$ and $q$ of the following form:
\begin{equation*}
\begin{cases} 
p=2^{\alpha} i, \quad q=2^{\beta}j+e \quad  \textrm{if p is even} \\
p=2^{\beta} j+e, \quad q=2^{\alpha}i \quad \textrm{if p is odd},
 \end{cases}
\end{equation*}
where $\alpha\geq1$, $\beta\geq2$, $e=\pm1$, $i$ and $j$ are odd natural numbers. He proved that the conjecture holds for these values of $p$ and $q$.\\

Miyazaki also obtained some lower and upper bounds for the solutions, in particular for $x\pm z$, by using $2$-adic  and p-adic valuation and also Baker theory. He showed the following theorem:
\begin{theorem}\label{prop:2.6}
If equation \eqref{eq:2.2} has a solution $(x,y,z)$ with even $x$ and $z$, then $x/2$ and $z/2$ must be odd. Also if $2\alpha=\beta+1$, then $x=y=z=2$.
\end{theorem}
He proved this result by using some published results, see \cite{Br}, \cite{Br2}, \cite{Ca2}, \cite{CD}, \cite{Co}, \cite{Da}, \cite{DM}.\\

In $2012$, Fujita and Miyazaki \cite{FuMi} proved the conjecture when $V\equiv 0\pmod{2^{r}}$ and $V\equiv \pm 2^{r}\pmod{U}$ for some non-negative integer $r$, and concluded also that the conjecture holds when $W\equiv -1\pmod{U}$.\\

Miyazaki \cite{Mi3} proved that the conjecture holds when $U\equiv \pm1\pmod{V}$ (see Theorem 1 of \cite{Mi3}), $W\equiv 1\pmod{V}$ (see Theorem 2 of \cite{Mi3}), and when $U-V= \pm1$ (as a consequence of these results).\\

Later, Fujita and Miyazaki \cite{FujMiy} proved another special case of the conjecture:
\begin{theorem}
If $k>1$ is a divisor of $V$ such that $k\equiv \pm 1 \pmod{U}$ and $V/k$ has no prime factor congruent to $1$ modulo $4$, then the conjecture holds.
\end{theorem}
They applied a lemma of Laurent \cite{Lau} which gives explicit lower bounds for a linear form in two logarithms to solve the difficulty when $k\equiv -1 \pmod{V}$ and $y=1$ in which case the obtained equation becomes $V=W^{z}-U^{x}$, which was considered by Pillai \cite{Pi}. They also showed that the conjecture holds when $p-q$ has a divisor congruent to $\pm 3 \pmod 8$ or $p+q$ has a divisor congruent to $5$ or $7 \pmod 8$.\\

Miyazaki, Yuan and Wu \cite{MiYuWu} established the conjecture when $V$ is even and either $U$ or $W$ is congruent to $\pm 1$ modulo the product of all prime factors of $V$.\\

In $2014$, Terai \cite{Ter4} completely settled the conjecture for $q=2$ without any assumptions on $p$ by using Theorem \ref{prop:2.6} and a lemma of Laurent \cite{Lau}.\\

In $2015$, Miyazaki and Terai \cite{MTer2} generalized equation \eqref{eq:2.2} by proving the following theorem:
\begin{theorem}
Let $q\equiv 2 \pmod{4}$ be a positive integer. Suppose that $q$ satisfies at least one of the following conditions:\\

$(i)$ $q/2$ is a power of an odd prime,\\

$(ii)$ $q/2$ has no prime factors congruent to $1$ modulo $8$,\\

$(iii)$ $q/2$ is a square.\\
Furthermore, suppose that $p>72q$. Then, equation \eqref{eq:2.2} has the unique solution $(x,y,z)=(2,2,2)$ in positive integers.
\end{theorem}

%#####################################################################################################%
\section{the Je\'{s}manowicz'  Conjecture for non-primitive Pythagorean Triples: $(kU)^{x}+(kV)^{y}=(kW)^{z}$}\label{sec:3}
%####################################################################################################%

Recall that on the statement of the Je\'{s}manowicz'  conjecture, there is a condition that the numbers $U,V,W$ satisfy the Pythagorean equation $U^{2}+V^{2}=W^{2}$. These triples were called primitive triples. When the gcd of $U,V$ and $W$ is greater than $1$, there was no study on corresponding version of Je\'{s}manowicz' conjecture until $1998$. Since then, several authors studied the more general equation
\begin{equation}\label{eq:3.1}
(kU)^{x}+(kV)^{y}=(kW)^{z}
\end{equation}
under several conditions with $k>1$ and $U^{2}+V^{2}=W^{2}$.\\

In $1998$, Deng and Cohen \cite{DeCo} obtained the following result:
\begin{theorem}
Let $U,V$ and $W$ be as in \eqref{eq:3.1}, $U$ is a prime power and $k$ is a positive integer such that $P(V)\mid k$ or $P(k) \nmid V$, where $P(k)$ is the product of all distinct prime divisors of $k$. Then, the only solution of \eqref{eq:3.1} is $x=y=z$.
\end{theorem}

In fact, they also completed the study for $U,V,W$ choosen as in \eqref{eq:2.1} for $1\leq n\leq 5$:
\begin{theorem}
For each of the Pythagorean triples $(U,V,W)=(3,4,5),(5,12,13)$, $ (7,24,25),(9,40,41), (11,60,61)$ and for any positive integer $k$, the only solution of the Diophantine equation \eqref{eq:3.1} is $x=y=z=2$.
\end{theorem}

Following Deng and Cohen's work, Le \cite{M1} gave the following more general result in $1999$:
\begin{theorem} \label{theo:3.3}
If $(x,y,z)$ is a solution of \eqref{eq:3.1} with $(x,y,z)\neq(2,2,2)$ , then one of the following conditions is satisfied:\\
$(i)$ $max\{x,y\}>min\{x,y\}>z$, $P(k)\mid W$, and $P(k)<P(W)$,\\
$(ii)$ $x>y>z$ and $P(k)\mid V$,\\
$(iii)$ $y>z>x$ and $P(k)\mid U$.
\end{theorem}
As a result of this, one can obtain the following corollaries:
\begin{corol}
If $(x,y,z)$ is a solution of \eqref{eq:3.1} with $(x,y,z)\neq(2,2,2)$, then $x,y$ and $z$  are distinct. 
\end{corol}
\begin{corol}
If  $P(k)$ does not divide any one of $U,V$ and $W$, then \eqref{eq:3.1} has only solution $(x,y,z)\neq(2,2,2)$. 
\end{corol}

$16$ years later in $2015$, Yang and Fu \cite{YF} simplified the conditions given in Theorem \ref{theo:3.3} by removing all conditions on $P(k)$. Meanwhile between $1999$ and $2015$, many mathematicians considered several specific cases of equation \eqref{eq:3.1}.\\

In $2013$, Yang and Tang \cite{YaTa} obtained some results when $U,V,W$ are related to Fermat primes:
\begin{theorem} 
Let $n$ be a positive integer. If $F_{n}=2^{2^{n}}+1$ is a Fermat prime, then for any positive integer $k$, the Diophantine equation
\begin{equation}\label{eq:3.2}
((F_{n}-2)k)^{x}+(2^{2^{n-1}+1}\, k))^{y}=(F_{n}k)^{z}
\end{equation}
has no solutions $(x,y,z)$ satisfying $z<min\{x,y\}$.
\end{theorem}
\begin{theorem} 
Let $n\geq 4$ be a positive integer and $F_{n}=2^{2^{n}}+1$. Then, for any positive integer $k$, \eqref{eq:3.2} has no solution other than $(x,y,z)=(2,2,2)$.
\end{theorem}

In $2014$, Tang and Weng \cite{Tang} generalized the above result and proved that the unique solution of \eqref{eq:3.2}, for any positive integers $n$ and $k$, is $(x,y,z)=(2,2,2)$.\\

The same year, Xinwen and Wenpeng \cite{XinWen} studied a special case of equation \eqref{eq:2.3}, where $U=2^{2m}-1$, $V=2^{m+1}$ and $W=2^{2m}+1$ with $p=2^{m}$ and $q=1$ in \eqref{eq:3.1}. They proved that the only solution of the equation
\begin{equation}\label{eq:3.3}
((2^{2m}-1)k)^{x}+(2^{m+1}k)^{y}=((2^{2m}+1)k)^{z},
\end{equation}
for any positive integers $m$ and $k$, is $(x,y,z)=(2,2,2)$. Actually, a special case of this equation, for $k=1$ and $m=\log_{2}2n$ was already completely solved by Lu \cite{Lu1} in $1959$. Also, the paper \cite{DeCo} by Deng and Cohen covered the case $m=1$ in equation \eqref{eq:3.3} and proved that the equation
\begin{equation*}
(3k)^{x}+(4k)^{y}=(5k)^{z}
\end{equation*}
has the unique solution $(x,y,z)=(2,2,2)$, for all $k>1$.\\

Finally, another special case for $m=2$ of equation \eqref{eq:3.3} was recently studied by Yang and Tang \cite{YaTa2}. They proved that the only solution of
\begin{equation*}
(15k)^{x}+(8k)^{y}=(17k)^{z}
\end{equation*}
is $(x,y,z)=(2,2,2)$, for $k\geq 1$.\\

In $2014$, Deng \cite{De1} considered another special case of equation \eqref{eq:3.3} by putting $m=s+1$ and accepting some divisibility conditions such as $P(U)\mid k$ or $P(k)\nmid U$, $s\geq 0$, and proved that $(x,y,z)=(2,2,2)$ is the only solution of the equation
\begin{equation} \label{eq:3.4}
((2^{2s+2}-1)k)^{x}+(2^{s+2}k)^{y}=((2^{2s+2}+1)k)^{z}.
\end{equation} 
This equation corresponds to the equation considered by Lu \cite{Lu1} with $p=2^{s}$, $s\geq 0$, $k=1$. In the same paper, Deng also omitted the divisibility conditions and proved that the only solution of \eqref{eq:3.4} is $(x,y,z)=(2,2,2)$, when $1\leq s\leq4$.\\

In $2015$, Ma and Wu \cite{Ma-Wu} considered the generalization of the equation of Lu \cite{Lu1}:
\begin{equation}\label{eq:3.5}
((4n^{2}-1)k)^{x}+(4kn)^{y}=((4n^{2}+1)k)^{z}.
\end{equation} 

If $P(4n^{2}-1)\mid k$, then they proved that the only solution of \eqref{eq:3.5} is $x=y=z=2.$ Another result they gave was on the case where $n=p^{\alpha}$, $\alpha\geq 0$, where $p$ is a prime congruent to $3$ modulo $4$. If $P(k)\nmid (4n^{2}-1)$, then they showed that the only solution is again $x=y=z=2$.\\

Very recently, Miyazaki \cite{Mi4} nicely proved Je\'{s}manowicz' Conjecture when $U$ or $V$ is a power of $2$, by extending the result of Tang and Weng \cite{Tang}.

%##############################################################%
\section{Some Bi-Products of The Je\'{s}manowicz' Conjecture}\label{sec:4}
%##############################################################%

While trying to prove or disprove the Je\'{s}manowicz' conjecture, several mathematicians considered different variants of it. The most important one of these variants is Terai's conjecture. The following slightly different version of Conjecture \ref{333} is known as the first Terai conjecture \cite{Ter94}:
\begin{conj}[\textbf{First Terai Conjecture}]\label{223}
If $a$, $b$, $c$, $p$, $q$, $r$ are fixed positive integers satisfying $a^{p}+b^{q}=c^{r}$ with $p,q,r\geq2$, then the Diophantine equation
\begin{equation*}
a^{x}+b^{y}=c^{z}
\end{equation*}
has only the positive integer solution $(x,y,z)=(p,q,r)$. 
\end{conj}
Several studies on this conjecture are as follows:\\

Terai \cite{Ter94}, Le \cite{Le97}, Cao and Dong \cite{CaDo1} considered the case $p=q=2,$ $r=3$ and 
$a=m^{3}-3m, \ b=3m^{2}-1, \ c=m^{2}+1,$ 
with $ 2\mid m $.\\

Terai \cite{Ter95}, Cao and Dong \cite{CaDo1} took the case $p=q=2$, $r=5$ and 
$$a = m( m^{4}-10m^{2}+5), \ b=5m^{4}-10m^{2}+1, \ c= m^{2}+1,$$ 
with $ 2\mid m$.\\

In $1999$, Cao \cite{Ca1} showed that first Terai conjecture is false. For example, from Nagell’s result \cite{Nagell}, we see that the equation $3^x + 2^y = 5^z$ has two solutions $(x,y,z)=(1,1,1), (2,4,2)$, and the equation $7^x +2^y = 3^z$ also has two solutions $(x, y, z) = (1, 1, 2), (2, 5, 4)$. Furthermore, if $a = 1$ or $b = 1$, then the conjecture is also false. And in his paper, Cao studied the case $p=q=2$, $2\nmid r$, $c\equiv 5 \ \pmod{8}$, $b\equiv 3 \pmod{4}$ and $ c $ is a prime power.\\ 

Same year, in same journal, Terai \cite {Ter2} modified his own conjecture as follows:

\begin{conj} [\textbf{Second Terai Conjecture}] \label{555} If $ a, \ b, \ c, \ p, \ q, \ r \in \mathbb{Z}^{+} $ are fixed integers satisfying $ a^{p}+b^{q}=c^{r} $, with $ p, \ q, \ r\geq2$ and $\left( a, \ b\right) =1$, then the Diophantine equation
\begin{equation}\label{eq:4.2}
a^{x}+b^{y}=c^{z}
\end{equation}
has only the positive integral solutions $\left( x,y,z \right) = \left( p,q,r \right)$, except for the following three cases $\left( with \ a \textless b  \right)$, where equation \eqref{eq:4.2} has only the following solutions, respectively: 
\[(a,b,c) = (2,3,5), \ (x,y,z) = (1,1,1), (4,2,2), \]
\[(a,b,c) = (2,3,7), \ (x,y,z) = (1,1,2), (5,2,4), \] 
\[(a,b,c) = (1,2,3), \ (x,y,z) = (m,1,1), (n,3,2), \]  
where $m$ and $n$ are orbitrary integers.   
\end{conj}
Before giving this conjecture, Terai considered the special case $ p=q=2 $ and $r$ an odd prime. \cite{Ter94},  \cite{Ter95}, \cite{Ter96}. \\

Terai also used a result on a lower bound for linear forms in two logarithms to the Diophantine equation
\begin{equation}\label{eq:4.3}
a^{x}+b^{n}=c^{z},
\end{equation} 
where $n$ is a given small positive integer. In fact, Terai applied a result of Laurent, Mignotte and Nesterenko \cite{LaMiNe} to obtain a lower bound for linear forms in two logarithms to deal with equation \eqref{eq:4.3}, where $n$ is a given ``relatively small"  positive integer. \\

In $2002$, Cao and Dong \cite{CD} considered the case $p=q=2$ and $ r\geq 3$ odd, $b\equiv 3 \pmod{4}$, $2 \parallel a$ and $b\geq25.1a$. They also proved that the equation $(2^{n}-1)^{x}+2^{y}=(2^{n}+1)^{z}$ has two solutions $(x,y,z)=(1,1,1)$ and $(2,n+2,2)$ for any $1<n\in\mathbb{Z^{+}}$. So they suggested that Terai's conjecture should be modified as follows: 
\begin{conj}[\textbf{Terai-Je\'{s}manowicz' Conjecture}]
If $a, \ b, \ c, \ p, \ q, \ r \in \mathbb{N} $ with $a^{p}+b^{q}=c^{r}, \ a, \ b, \ c, \ p, \ q,\ r\geq2$ and $(a,b)=1$, then Diophantine equation \eqref{eq:4.2} has only the solution $ \left( x,y,z\right) =\left( p,q,r\right)$, with $ x,y,z>1 $.	 
\end{conj}
Note that this conjecture which is called as Terai-Je\'{s}manowicz' conjecture by the authors is yet another generalization of the Je\'{s}manowicz' conjecture obtained when $ p=q=r=2 $.

Same year, the same authors together with Li \cite{CaDoLi} proved the following result which is a slight generalization of the above result by using a result of Bilu, Hanrot and Voutier \cite{BHV}:
\begin{theorem}
With the above notation and $b$ is an odd prime power, then equation \eqref{eq:4.2} has the unique solution $(x,m,n)=(a,2,r)$.
\end{theorem}

Next year, Le \cite{Le4} showed that "Terai-Je\'{s}manowicz' conjecture" is also false. For example, if $a=2$, $b=2^n-1$, $c=2^n+1$, where $n$ is a positive integer with $n>2$, then $a, b, c$ satisfy $max(a,b,c)>7$ and $a^{n+2}+b^2=c^2$, but the equation $a^x+b^y=c^z$ has two solutions $(x,y,z)=(1,1,1)$ and $(n+2,2,2)$ where $a, b, c$ are fixed coprime positive integers with $min(a,b,c)>1$ and $x,y,z$ are integers. This implies that there exist infinitely many counterexamples to Terai-Je\'{s}manowicz' conjecture. Le suggested the following conjecture:
\begin{conj}
The equation $a^x+b^y=c^z$ has at most one solution $(x,y,z)$ with $min(x,y,z)>1$.
\end{conj}
The above mentioned conjecture was first proposed by Le \cite{Le5} for primes $a$, $b$ and $c$. It was proved for some special cases. But, in general, the problem has not been solved yet. In  \cite{Le4}, Le considered the case $a^{2}+b^{2}=c^{r}$, $gcd(a,b)=1$, $2\nmid a$, $2\mid b$, $r>1$, $2\nmid r$ where $a, b, c$ positive integers and proved if $a>b$, $a\equiv 3 \pmod{4}$, $b\equiv 2 \pmod{4}$ and $a/b>(e^{r/1856}-1)^{-1/2}$, then the equation $a^x+b^y=c^z$ has only the solution $(x,y,z)=(2,2,r).$\\

In $2006$, same author \cite{Le1} changed the condition $a^{2}+b^{2}=c^{2}$ to $a^{2}+b^{2}=c^{r}$, $r$ odd and proved the following result:
\begin{theorem}\label{Theorem 4.2} 
Let $a,b,c,r\in \mathbb{Z^{+}}$ such that $a^{2}+b^{2}=c^{r}$, $min(a,b,c,r)>1$, $gcd(a,b)=1$, $a$ is even and $r$ is odd. If $b\equiv 3 \pmod{4}$ and either $b$ or $c$ is an odd prime power, then the equation \eqref{eq:4.2} has the unique solution $(x,m,n)=(a,2,r)$ with $min(m,n)>1$.
\end{theorem}     

In $2009$, Cipu and Mignotte \cite{CiMi} improved the conditions given in the papers \cite{Ca1}, \cite{CD}, \cite{Le2}, and \cite{Ter2}. Moreover, they proved three main theorems under the conditions $a\equiv 2 \pmod{4}$, $b\equiv 3 \pmod{4}$, $(a,b)=1$, $r>1$ odd and $a^{2}+b^{2}=c^{r}$. They also proved that if there was a counterexample to Terai's conjecture, then each of $a$, $b$ and $c$ would have at least two prime divisors.\\

In $2011$, Miyazaki \cite{Mi8} considered the case $r>2$ is even and proved Terai-Je\'{s}manowicz' conjecture without any assumptions on $a$, $b$ and $c$.\\

The reader can also consult \cite{YHT-2012} for the progress made in the tentatives to solve the equation. Moreover, Yang, He, and Togb\'e improved the bounds obtained by the previous authors under some conditions.\\
 
In $2010$, Miyazaki \cite{Mi7} stated a different version of Je\'{s}manowicz' conjecture, which he called the shuffle variant of Je\'{s}manowicz' conjecture:\\
\begin{conj}[\textbf{Shuffle Variant of Je\'{s}manowicz' Conjecture}]\label{Conj 4.4}
Let (a,b,c) be a primitive Pythagorean triple with $b$ even. If $c=b+1$, then the equation
\begin{equation} \label{eq:4.4}
c^{x}+b^{m}=a^{n}
\end{equation}has the unique solution $(x,y,z)=(1,1,2)$. If $c>b+1$, then  there are no solutions.  
\end{conj} 
Moreover, he proved that if $c\equiv 1 \pmod{b}$, then Conjecture \ref{Conj 4.4} is true.\\
 
In \cite{Mi3}, he studied the same equation in Theorem \ref{Theorem 4.2} with one change that $r$ is even. Moreover, he similarly proposed a variant of Terai's conjecture, which he parallelly called as the shuffle variant of Terai's conjecture:
\begin{conj}[\textbf{Shuffle Variant of Terai's Conjecture}]
Let $p,q,r \geq 2$ be integers and let $a,b,c$ be pairwise relatively prime positive integers such that $ a^{p}+b^{q}=c^{r}$. If  $\left( a,b,c\right)$ is different from $\left( 2,7,3\right)$ and $\left( 2,2^{p-2}-1,2^{p-2}+1\right)$ where $p\geq 3$, then equation \eqref{eq:4.4} has the unique solution $\left(x,m,n\right) =\left( 1,1,p\right)$ if $q=r=2$ and $c=b+1$ and no solution otherwise. 	  
\end{conj} 
He gave some conditions for which this conjecture is true and also proved the first part of it. \\

In $2014$, Le, Togb\'e and Zhu \cite{LTZ} considered the equation \eqref{eq:4.2} and proved that when $m > \max\{10^{15},2r^{3}\}$, the equation \eqref{eq:4.2} has the unique solution $\left( x,y,z\right) =\left( 2,2,r\right)$. The same year, Terai \cite{Ter6} considered the equation
\begin{equation}\label{eq:4.5}
x^{2}+q^{m}=c^{n}, 
\end{equation} where $q^{t}+1=2c^{s}$ with $q$ prime and $s=1,2$. He proved three theorems on the solutions of \eqref{eq:4.5} under some conditions.

%################################%
\section{Several Variants}\label{sec:5}
%################################%

In the last $15$ years, there has been a series of papers dealing  with particular classes of numbers in place of $a$, $b$, and/or $c$.
In $2001$, Terai and Takakuwa \cite{TerTa} called the positive integers $a,b,c$ such that $a^{2}+ ab+b^{2}=c^{2}$ as \textbf{Eisenstein numbers} and taking this new condition instead of $a^{2}+b^{2}=c^{2}$. They studied the equation
 \begin{equation}\label{eq:5.1} 
a^{2x}+a^{x}b^{y}+b^{2y}=c^{z}
\end{equation}
and conjectured that $\left( x,y,z\right) =\left( 1,1,2\right)$ is the unique solution when $(a,b)=1$. They also showed that when $a$ or $b$ is a prime power, their conjecture holds. \\

In $2009$, He and Togb\'e \cite{HeTo} studied the exponential Diophantine equation
\begin{equation}\label{eq:5.2}
n^{x}+(n+1)^{y}=(n+2)^{z}
\end{equation} 
and showed that all solutions are $(n,x,y,z)=(1,t,1,1),(1,t,3,2),(3,2,2,2)$. In fact, equation \eqref{eq:5.2} is a generalization of the equation $3^x+4^y=5^z$ studied by Sierpi\'nski~\cite{Si}. The general equation was first studied by Leszczy\'nski \cite{Leszczynski:1959} and Makowski \cite{Makowski:1966}.  Makowski extended Leszczy\'nski's work and found the solutions when $y=1$ for $1 \leq n \leq 48$.
The equation was solved when $xyz=0$ (see Theorem 2, \cite{Makowski:1966}).\\
 
In $2012$, Terai \cite{Ter2012 } studied the equation
\begin{equation}\label{eq:5.3}
(4m^{2}+1)^{x}+(5m^{2}-1)^{y}=(3m)^{z}
\end{equation} and by means of Baker's theory, he showed that if $m\leq 20$ and
$m\not\equiv 3 \pmod{6}$, then the only solution is $(x,y,z)=(1,1,2)$. In $2014$, Su and Li \cite{SuLi} generalized this result to $ m>90 $ and $ m\equiv 3 \pmod{6}$.\\

Recently, Bert\'{o}k \cite{Ber} completely solved equation \eqref{eq:5.3} by using exponential congruences to deal with the remaining values, i.e.  $20<m\leq 90$ and $m\equiv 3 \pmod{6}$.\\
 
In $2012$, Miyazaki and Togb\'{e} \cite{MT} studied another version of this conjecture by considering the equation
\begin{equation}
(2am-1)^{x}+(2m)^{y}= (2am+1)^{z},
\end{equation}   
for any fixed positive integer $a>1$, and proved that all solutions are $ (m,x,y,z)=(2a,2,2,2), (1,1,1,1),$ giving also some further results for particular values of $a$. By means of these results, they showed that the equation 
\begin{equation}\label{eq:5.4}
b^{x}+2^{y}=(b+2)^{z}
\end{equation}
has only the solution is $(x,y,z) = (1,1,1)$, with only one exception, when $b=89$, it has also $(1,13,2)$ as a solution.\\   

In $2013$, Tang and Yang \cite{YaTa3} generalized equation \eqref{eq:5.4} to
\begin{equation} \label{eq:5.5}
(bn)^{x}+(2n)^{y}=((b+2)n)^{z},
\end{equation} 
with $ b\geq 5 $ is an odd integer and showed that if $(x,y,z)\neq (1,1,1)$ is a solution, then $y<z<x$ or $x\leq z<y$.\\

In $2014$, Miyazaki and Terai \cite{MTer} studied the equation
\begin{equation} \label{eq:5.6}
(m^{2}+1)^{x}+(cm^{2}-1)^{y}=(am)^{z},
\end{equation} 
when $a,c,m$ are positive integers with $ a\equiv 3,5 \pmod{8}$. They showed that if $ 1+c=a^{2}$ and $m\equiv \mp 1 \pmod{a}$, then  equation \eqref{eq:5.6} has the unique solution $(x,y,z)=(1,1,2)$, except for the case $(m,a,c)=(1,3,8)$ where $(5,2,4)$ is also a solution. \\

In $2015$, Miyazaki \cite{Mi5} answered a question proposed by Terai, \cite{Ter5} during ICM $2002$, which states that equation \eqref{eq:2.2}, with $U=F_{n}$, $V=F_{n+1}$ and $W=F_{2n+1}$ Fibonacci numbers, has the unique solution $(x,y,z)=(2,2,1)$ in positive integers. Moreover, he showed another result as follows: 
\begin{theorem}
For each $n\geq 3$, the exponential Diophantine equation
\begin{equation}
(F_{n})^{x}+(F_{2n+2})^{y}=(F_{n+2})^{z}
\end{equation}
has the unique solution $(x,y,z)=(2,1,2)$ in positive integers.
\end{theorem} 
    
In $2015$, Terai and Hibino \cite{TerHib} considered another exponential Diophantine equation
\begin{equation} \label{eq:5.7}
(12m^{2}+1)^{x}+(13m^{2}-1)^{y}=(5m)^{z}
\end{equation} and using Baker's theory they proved that it has the unique positive integer solution $(x,y,z)=(1,1,2)$ when $m\not\equiv 17, 33 \pmod{40}$.\\

Very recently, Miyazaki, Togb\'{e} and Yuan \cite{Mi6} generalized \eqref{eq:5.4} by taking $c$ instead of $2$ and determined all positive integer solutions of
\begin{equation} \label{eq:5.8}
b^{x}+c^{y}=(b+2)^{z},
\end{equation}
for any coprime positive integers $b$, $c$ with $b\equiv -1 \pmod{c}$ by means of Baker's theory.\\
     
Terai and Hibino \cite{TerHib2} gave the most recent results by considering the following exponential Diophantine equation
\begin{equation}\label{eq:5.9} 
(3pm^{2}-1)^{x}+(p(p-3)m^{2}+1)^{y}=(pm)^{z},
\end{equation}
where $ p\equiv 1 \pmod{4}$ is a prime and $m$ is a positive integer. They proved the following result.\\
\begin{theorem}
Let $m$ be a positive integer with $m \not\equiv 0 \ (mod \ 3) $. Let $p$ be a prime with $p\equiv 1 \pmod{4}$. Moreover, suppose that if $m\equiv 1 \pmod{4}$, then $p<3784$. Then, equation \eqref{eq:5.9} has the unique positive integer solution $(1,1,2)$. 
\end{theorem}
They noticed that for $p=5$, it is possible to solve this equation without any  assumption. \\

One must remark that \eqref{eq:5.9} is a particular member of a more general class of exponential Diophantine equations 
\begin{equation} \label{eq:5.10}
t^{x}+(k^{2}-t)^{y}=k^{z},
\end{equation} with $k=pm$ and $t=3pm^{2}-1$. 

%##############################################################%
\section{An "elementary" method to prove a new result}\label{sec:6}
%##############################################################%

In this section, for the sake of brevity, the authors consider \eqref{eq:3.1} with $(U,V,W)=(20,99,101)$ and they solve the Diophantine equation
\begin{equation}\label{eq:2}
(20k)^x+(99k)^y=(101k)^z.
\end{equation} 

First, we briefly recall the results that led the authors to deal with this case.
Recently, several authors showed that Je\'{s}manowicz' conjecture is true for $2\leq n\leq4$ and $n=8$ where $(U,V,W)=(4n,4n^{2}-1,4n^{2}+1)$ for equation \eqref{eq:3.1}. One can have a look at \cite{De1}, \cite{Tang}, \cite{YaTa2}, and \cite{Zhiju}. Recall also that the very first case where $n=1$ was settled in 1959 by Lu \cite{Lu1}.

Therefore, the natural next step is to consider the case $n=5$. One obtains equation \eqref{eq:2} in this case. Here is our main result:

\begin{theorem}\label{thm:main}
Let $k$ be any positive integer. Then, Diophantine equation \eqref{eq:2} has only the solution $(x,y,z)=(2,2,2)$.
\end{theorem}

We organize the proof of our main result as follows. In Subsection~\ref{subsec:6.1}, we recall two useful results that we will use for the proof of Theorem \ref{thm:main}. Subsection~\ref{subsec:6.2} is devoted for this proof. We know that any solution $(x, y, z)\neq (2, 2, 2)$ of equation \eqref{eq:2} verifies $z< \max\lbrace x,y\rbrace.$ Therefore, we will consider two cases: $x<y$ and $x>y$. For the sake of completeness, we will give all details for each the above cases.
\vspace{5mm}

%%%%%%%%%%%%%%%%%%%%%%%%%%
\subsection{Lemmas}\label{subsec:6.1}
%%%%%%%%%%%%%%%%%%%%%%%%%%

In this subsection, we will only recall the two useful results necessary for the proof of our main result. 

\begin{lemma} \label{lem:1} (\cite{Lu1}) 
The Diophantine equation $$(4n^2-1)^x+(4n)^y=(4n^2+1)^z$$ has only the solution $(x,y,z)=(2,2,2)$. 
\end{lemma}
\begin{lemma} \label{lem:2}(\cite{DeCo}, Lemma 2) 
If $z\ge max\lbrace x,y\rbrace$, then the Diophantine equation $U^x+V^y=W^z$, where $U, V$, and $W$ are any positive integers (not necessarily relative prime) such that $U^2+V^2=W^2$, has no solution other than $(x,y,z)=(2,2,2)$. 
\end{lemma}

\vspace{5mm}

%%%%%%%%%%%%%%%%%%%%%%%%%%%%%
\subsection{Proof of Theorem $6.1$}\label{subsec:6.2}
%%%%%%%%%%%%%%%%%%%%%%%%%%%%%

Recall that according to Lemma \ref{lem:1}, the only solution $(x,y,z)$ of 
$$20^x+99^y=101^z$$ 
is $(2,2,2).$ On the contrary, we assume that equation \eqref{eq:2} has at least another solution than $(x, y, z)\neq (2, 2, 2)$. Using Lemma \ref{lem:2}, we conclude that $n\geq2$ and $z< max\lbrace x,y\rbrace.$ Moreover, from the result in \cite{M1}, we won't consider the cases $x=y$, $x=z$, and $y=z$.
\\
\\
%***************************%
\textbf{Case 1. $x<y.$}
%***************************%
Because of the condition $z< max\lbrace x,y\rbrace$, we will consider two subcases.
\\
\\
\textbf{Subcase 1.1 $z<x<y.$} Then,
\begin{equation} \label{eq:3}
k^{x-z}(20^x+99^yk^{y-x})=101^z.
\end{equation}  

If $(k,101)=1$, then equation \eqref{eq:3} and $k\geq2$ imply $x=z$, which is a contradiction.\\

If $(k,101)=101$, then we write $k=101^{s} n_{1}$, where $s\geq1$ and $101\nmid n_{1}.$ Using \eqref{eq:3}, we have
\begin{equation} \label{eq:4}
101^{s(x-z)}n_{1}^{x-z}(20^x+99^y101^{s(y-x)}n_{1}^{y-x})=101^z.
\end{equation}  
It follows that $n_{1}^{x-z}\mid 101^{z}$, thus $n_{1}=1.$ One can see that equation \eqref{eq:4} implies
\begin{equation*}
101^{s(x-z)}(20^x+99^y101^{s(y-x)})=101^z.
\end{equation*}
Then $s(x-z)<z$ and $20^{x}+99^y101^{s(y-x)}=101^{z-s(x-z)}$, thus $101\mid 20^{x}$. We deduce a contradiction.  
\\
\\
\textbf{Subcase 1.2 $x < z<y.$} Then, we get
\begin{equation} \label{eq:5}
20^{x}+99^{y}k^{y-x}=101^zk^{z-x}.
\end{equation} 

If $(k,20)=1$, thus equation \eqref{eq:5} and $k>2$ imply that $x=z<y.$ We deduce a contradiction to the fact that $x<z$. Therefore, we suppose $(k,20)>1$. We write $k=2^{r}5^{s}n_{1}$ where $r+s\geq1$ and $(10,n_{1})=1.$ We transform equation  \eqref{eq:5} to obtain 
\begin{equation} \label{eq:7}
20^{x}=2^{r(z-x)}5^{s(z-x)}n_{1}^{z-x}[101^{z}-99^{y}2^{r(y-z)}5^{s(y-z)}n_{1}^{y-z}].
\end{equation}
\\

(i) If $r\geq1$, $s=0$, then write $k=2^{r}n_{1}$, where $(10,n_{1})=1$. Equation \eqref{eq:7} becomes 
\begin{equation*}
20^{x}=2^{r(z-x)}n_{1}^{z-x}[101^{z}-99^{y}2^{r(y-z)}n_{1}^{y-z}].
\end{equation*} 
We deduce that $r(z-x)=2x$ and
\begin{equation*}
5^{x}=n_{1}^{z-x}[101^{z}-99^{y}2^{r(y-z)}n_{1}^{y-z}].
\end{equation*} 
As $(n_{1},10)=1$, then $n_{1}=1$ and
\begin{equation}\label{eq:8}
101^{z}-5^{x}=99^{y}2^{r(y-z)}.
\end{equation}
Thus, by consideration modulo $33$, one can see that $2^{z}\equiv 5^{x}\pmod{33}$ has only the solutions $(z,x)=(8,2),(18,2),(28,2).$ 
Put $x=2$ and $z=2z_{1},$ $z_{1} >0.$ By \eqref{eq:8}, we have  
\begin{equation*}
(101^{z_{1}}-5)(101^{z_{1}}+5)=99^{y}2^{r(y-z)}.
\end{equation*}
We know that that $(101^{z_{1}}-5,101^{z_{1}}+5)=2$, then
\begin{equation}\label{eq:9}
11^{y}\mid (101^{z_{1}}+5)
\end{equation} 
or
\begin{equation}\label{eq:10}
11^{y}\mid (101^{z_{1}}-5).
\end{equation}
However,
\begin{equation*}
11^{y}>11^{z}>(101+5)^{z_{1}}>101^{z_{1}}+5>101^{z_{1}}-5,
\end{equation*} 
which contradicts \eqref{eq:9} and \eqref{eq:10}. 
\\

(ii) If $r=0$, $s\geq1$, then we put $k=5^{s}n_{1}$, where $(10,n_{1})=1.$ Using \eqref{eq:7}, we get 
\begin{equation*}
20^{x}=5^{s(z-x)}n_{1}^{z-x}[101^{z}-99^{y}5^{s(y-z)}n_{1}^{y-z}].
\end{equation*} 
Then, $s(z-x)=x.$ We obtain
\begin{equation*}
2^{2x}=n_{1}^{z-x}[101^{z}-99^{y}5^{s(y-z)}n_{1}^{y-z}].
\end{equation*} 
As $(n_{1},10)=1$, we see that $n_{1}=1$ and
\begin{equation}\label{eq:11}
101^{z}-2^{2x}=99^{y}5^{s(y-z)}.
\end{equation}
The congruence modulo $11$ of the above equation gives $2^{z}\equiv 4^{x}\pmod{11}$. So $z\equiv 0\pmod{2}.$ We write $z=2z_{1}, z_{1} >0 $ and equation \eqref{eq:11} becomes   
\begin{equation*}
(101^{z_{1}}-2^{x})(101^{z_{1}}+2^{x})=99^{y}5^{s(y-z)}.
\end{equation*}
As $(101^{z_{1}}-2^{x},101^{z_{1}}+2^{x})=1$, we get
\begin{equation}\label{eq:12}
11^{y}\mid (101^{z_{1}}-2^{x})
\end{equation} 
or
\begin{equation}\label{eq:13}
11^{y}\mid (101^{z_{1}}+2^{x}).
\end{equation}
However, the following inequalities
\begin{equation}
11^{y}>11^{z}>(101+4)^{z_{1}}>101^{z_{1}}+2^{2z_{1}}>101^{z_{1}}+2^{x}>101^{z_{1}}-2^{x}
\end{equation} 
contradict \eqref{eq:12} and \eqref{eq:13}.
\\

(iii) If $r\geq1$, $s\geq1.$ Thus, $r(z-x)=2x,$ $x=s(z-x).$ From equation \eqref{eq:7}, one deduces that
\begin{equation*}
n_{1}^{z-x}[101^{z}-99^{y}2^{r(y-z)}5^{s(y-x)}n_{1}^{y-z}]=1.
\end{equation*} 
If $x=z$, then clearly $x=0$, which is a contradiction. Therefore, $x<z$ and $n_{1}=1.$ We obtain
\begin{equation}\label{eq:14}
99^{y}2^{r(y-z)}5^{s(y-x)}=101^{z}-1.
\end{equation}
The congruence modulo $3$ of the above equation gives $2^{z}\equiv 1\pmod{3}$ and then $z\equiv 0\pmod{2}.$ By consideration of equation \eqref{eq:14} modulo $17$, we have $101^{z}-1\equiv (-1)^{z}-1\equiv 0 \pmod{17}.$ Then, $17\mid (101^{z}-1)$. But it is clear that $17\nmid 99^{y}2^{r(y-z)}5^{s(y-x)}$, which is a contradiction.  This completes the proof for the first case.
\\
\\
%************************************%
\textbf{Case 2 $x>y.$} 
%************************************%
Again here, with the condition $z< \max\lbrace x,y\rbrace$, we have to consider two subcases. 
\\
\\
%&&&&&&&&&&&&&&&&&&&&&&&%
\textbf{Subcase 2.1 $z < y<x.$} 
%&&&&&&&&&&&&&&&&&&&&&&&%
Then, we have
\begin{equation} \label{eq:17}
k^{y-z}(20^xk^{x-y}+99^y)=101^z.
\end{equation}  

If $(k,101)=1$, then using \eqref{eq:17} and $k\geq2$, we obtain $y=z$. This contradicts the fact that $z<y$. 
If $(k,101)=101$, then put $k=101^{r}n_{1},$ where $r\geq1$ and $101\nmid n_{1}.$ Equation \eqref{eq:17} becomes
\begin{equation*}
n_{1}^{y-z}101^{r(y-z)}(20^{x}n_{1}^{x-y}101^{r(x-y)}+99^{y})=101^{z}.
\end{equation*} 
Because $(n_{1},101)=1$ and $(20^{x}n_{1}^{x-y}101^{r(x-y)}+99^{y},101)=1,$ we see that $r(y-z)=z$ and then $n_{1}^{y-z}(20^{x}n_{1}^{x-y}101^{r(x-y)}+99^{y})=1.$ This is also impossible.
\\
\\
%&&&&&&&&&&&&&&&&&&&&&&&%
\textbf{Subcase 2.2  $y<z<x.$} 
%&&&&&&&&&&&&&&&&&&&&&&&%
Then, equation \eqref{eq:2} is transformed into
\begin{equation} \label{eq:21}
99^{y}=k^{z-y}(101^z-20^{x}k^{x-z}).
\end{equation} 

If $(k,99)=1$, then by equation \eqref{eq:21} and $k>2$, we have $y=z$. This is a contradiction.
If $(k,99)>1$, then we write $k=3^{r}11^{q}n_{1}$, where $r+q\geq1$ and $(99,n_{1})=1.$ Equation \eqref{eq:21} implies  
\begin{equation} \label{eq:22}
99^{y}=3^{r(z-y)}11^{q(z-y)}n_{1}^{z-y}[101^{z}-20^{x}3^{r(x-z)}11^{q(x-z)}n_{1}^{x-z}].
\end{equation}
We will study all possibilities.
\\

(i) If $r\geq1$, $q=0$, then $k=3^{r}n_{1}$ and equation \eqref{eq:22} becomes
\begin{equation*}
99^{y}=3^{r(z-y)}n_{1}^{z-y}[101^{z}-20^{x}3^{r(x-z)}n_{1}^{x-z}].
\end{equation*} 
Then, $r(z-y)=2y$ and $(n_{1}, 99)=1$ imply $n_{1}=1.$ We get
\begin{equation}\label{eq:23}
20^{x}3^{2y}=101^{z}-11^{y}.
\end{equation} 
By consideration modulo $4$, we have $(-1)^{y}\equiv 1\pmod{4}$ and then $y$ is even, i.e. $y=2y_{1},y_{1} > 0.$ Now modulo $6$, we get $(-1)^{z}\equiv (-1)^{2y_{1}}\pmod{6}$. So $z$ must be even. Put $z=2z_{1},z_{1} > 0.$  From equation \eqref{eq:23}, we obtain
\begin{equation}\label{eq:24}
20^x3^{r(x-2z_1)}=(101^{z_{1}}-11^{y_{1}})(101^{z_{1}}+11^{y_{1}}).
\end{equation} 
As $(101^{z_{1}}-11^{y_{1}},101^{z_{1}}+11^{y_{1}})=2,$ we conclude that $101^{z_{1}}-11^{y_{1}}\equiv 0\pmod{5}.$ Hence, we have two possibilities: 
\begin{equation}\label{eq:25}
5^{x}2^{2x-1}\mid 101^{z_{1}}-11^{y_{1}} \; \mbox{ and }\;  2\mid 101^{z_{1}}+11^{y_{1}}
\end{equation} 
or
\begin{equation}\label{eq:26}
5^{x}2\mid 101^{z_{1}}-11^{y_{1}} \; \mbox{ and }\; 2^{2x-1}\mid 101^{z_{1}}+11^{y_{1}}.
\end{equation}
However, inequalities
$$2^{2x-1}5^{x}>2^{2z-1}5^{z}>2^{3z_{1}}5^{2z_{1}}>(101+11)^{z_{1}}>101^{z_{1}}+11^{y_{1}}>101^{z_{1}}-11^{y_{1}}$$
contradict \eqref{eq:25}. Therefore, we use \eqref{eq:26} to see that $101^{z_{1}}-11^{y_{1}}\equiv 1-(-1)^{y_{1}}\equiv 2\pmod{4}$. So $y_{1}$ is odd. If $3^{r(x-z)}\mid 101^{z_{1}}-11^{y_{1}}$, then from \eqref{eq:24}, we have
\begin{equation*}
101^{z_{1}}-11^{y_{1}}=2\cdot 5^{x}3^{r(x-z)},\quad 101^{z_{1}}+11^{y_{1}}=2^{2x-1}. 
\end{equation*} 
Hence, we get $11^{y_{1}}=2^{2x-2}-5^{x}3^{r(x-z)}.$ By consideration modulo $3$, we have $(-1)^{y_{1}} \equiv 1\pmod{3}$. This means that $y_{1}$ is even. It contradicts the fact that $y_{1}$ is odd. So $3^{r(x-z)}\mid 101^{z_{1}}+11^{y_{1}}.$ By \eqref{eq:24} and \eqref{eq:26}, we have
\begin{equation*}
101^{z_{1}}-11^{y_{1}}=2\cdot 5^{x}, 101^{z_{1}}+11^{y_{1}}=2^{2x-1}3^{r(x-z)}.
\end{equation*} Then we get
\begin{equation}\label{eq:27}
101^{z_{1}}=5^{x}+4^{x-1}3^{r(x-z)}
\end{equation} and
\begin{equation}\label{eq:28}
11^{y_{1}}=4^{x-1}3^{r(x-z)}-5^{x}.
\end{equation} 
If we consider equation \eqref{eq:27} modulo $12$, then we have $5^{z_{1}}\equiv 5^{x}\pmod{12}$, i.e. $5^{z_{1}-x}\equiv1\pmod{12}$. Then, we have two possibilities:\\

$A)$ {\bf Both $z_{1}, x$ are even}. Then, from \eqref{eq:27} and \eqref{eq:28}, we get
\begin{equation}\label{eq:35}
101^{z_{1}}-11^{y_{1}}=2\cdot 5^{x}. 
\end{equation} 
By reducing the above equation modulo $11,$, we get
\begin{equation*}
2^{z_{1}-1}\equiv 5^{x}\pmod{11}.
\end{equation*} 
As $5^{x}$ takes values $1, 3, 4, 5$, and $9$ modulo $11$, $2^{z_{1}-1}$ must also have these values modulo $11.$ In all cases $z_{1}$ must be odd, which contradicts the assumption that $z_{1}$ is even.\\

$B)$ {\bf Both $z_{1}, x$ are odd}. Then, reducing equation \eqref{eq:35} modulo $12$, we get a contradiction as $(-1)^{y_{1}}\equiv 7 \pmod{12}$.
\\

(ii) If $r=0$, $q\geq1$, then $k=11^{q}n_{1}.$ Therefore, equation \eqref{eq:22} becomes
\begin{equation*}
99^{y}=11^{q(z-y)}n_{1}^{z-y}[101^{z}-20^{x}11^{q(x-z)}n_{1}^{x-z}].
\end{equation*} 
Then $q(z-y)=y.$ As $(n_{1},99)=1,$ we get $k=1$ and then we have
\begin{equation}\label{eq:29}
20^{x}11^{q(x-z)}=101^{z}-3^{2y}.
\end{equation} 
Considering \eqref{eq:29} modulo $11$, we obtain $2^{z}\equiv(-2)^{y}\pmod{11}.$ We deduce that both $y, z$ are even. Put $z=2z_{1}, z_1 >0.$ We use equation \eqref{eq:29} to obtain  
\begin{equation}\label{eq:30}
20^{x}11^{q(x-z)}=(101^{z_{1}}-3^{y})(101^{z_{1}}+3^{y}).
\end{equation} 
As $(101^{z_{1}}-3^{y},101^{z_{1}}+3^{y})=2,$ we conclude that $101^{z_{1}}-3^{y}\equiv0\pmod{5}$, where $y\equiv0\pmod{4}.$ Put $y=4y_{1}$. Hence, we have two possibilities:
\begin{equation}\label{eq:31}
5^{x}2^{2x-1}\mid (101^{z_{1}}-3^{y}) \; \mbox{ and }\; 2\mid (101^{z_{1}}+3^{y})
\end{equation} or
\begin{equation}\label{eq:32}
5^{x}2\mid (101^{z_{1}}-3^{y}) \; \mbox{ and }\;  2^{2x-1}\mid (101^{z_{1}}+3^{y}).
\end{equation}
However, inequalities $2^{2x-1}5^{x}>2^{2z-1}5^{z}>2^{3z_{1}}5^{2z_{1}}>(101+9)^{z_{1}}>101^{z_{1}}+9^{z_{1}}=101^{z_{1}}+3^{2z_{1}}>101^{z_{1}}+3^{y}>101^{z_{1}}-3^{y}$ contradict the first relation of \eqref{eq:31}. Now, if $11^{q(x-z)} \mid (101^{z_{1}}+3^{y}),$ then from \eqref{eq:30} and \eqref{eq:32}, we deduce that 
\begin{equation*}
101^{z_{1}}+3^{y}=2^{2x-1}11^{q(x-z)}
\end{equation*}
and 
\begin{equation*}
101^{z_{1}}-3^{y}=2\cdot5^{x}.
\end{equation*}
Taking the difference of the above equations gives 
\begin{equation*}
3^{y}=2^{2x-2}11^{q(x-z)}-5^{x}.
\end{equation*}
By consideration modulo $4$, we have $3^{y}\equiv-1\pmod{4}$ So $y$ is odd. This contradicts the fact that $y$ is even. Therefore, $11^{q(x-z)} \mid (101^{z_{1}}-3^{y})$ and from \eqref{eq:30}, we have     
 \begin{equation*}
101^{z_{1}}-3^{y}=2\cdot5^{x}11^{q(x-2z_{1})} 
 \end{equation*}
and 
\begin{equation*}
101^{z_{1}}+3^{y}=2^{2x-1}.
\end{equation*}
We deduce 
\begin{equation}\label{eq:33}
101^{z_{1}}=2^{2x-2}+5^{x}11^{q(x-2z_{1})}
\end{equation}
and 
\begin{equation}\label{eq:34}
3^{y}=2^{2x-2}-5^{x}11^{q(x-2z_{1})}.
\end{equation}
Using \eqref{eq:33}, we have $2^{z_{1}}\equiv2^{2x-2}\pmod{11}$, which implies that $z_{1}$ is even. From \eqref{eq:34}, we get $(-1)^{4y_{1}}\equiv-(-1)^{q(x-2z_{1})}\pmod{4}$, so both $q$ and $x$ are odd. Finally, considering \eqref{eq:33} modulo $4$, we get $1\equiv (-1)^{qr} \equiv -1 \pmod{4}$ as $q, x$ are odd. This is a contradiction.
\\

(iii) If $r\geq1, q\geq1$ then from equation \eqref{eq:22}, we obtain $r=2q$ as $r(z-y)=2q(z-y)$. Hence,  
\begin{equation*}
1=n_{1}^{z-y}[101^{z}-20^{x}3^{r(x-z)}11^{q(x-z)}n_{1}^{x-z}].
\end{equation*}
As $y<z,$ we have $n_{1}=1$ and 
$$1=101^{z}-20^{x}3^{r(x-z)}11^{q(x-z)}.$$
Considering the above equation modulo $3$, we see that $1\equiv2^{z}\pmod{3}$, which implies that $z$ is even. Now, $101^{z}-1\equiv(-1)^{z}-1\equiv0\pmod{17}.$ Hence, $17\mid(101^{z}-1)$. But $17\nmid20^{x}3^{r(x-z)}11^{q(x-z)}$. This is a contradiction and completes the proof of Theorem \ref{thm:main}.

\subsection*{Acknowledgements}
The authors would  like to thank Professors Min Tang and Huilin Zhu for providing them some references. The first and third authors were supported by the research fund of Uludag University project no: F-2015/23, F-2016/9. The second author was supported by the research fund of Uludag University project no: F-2015/18. The fourth author thanks Purdue University Northwest for the support.

\end{document}